\documentclass[11pt,reqno]{amsart}

\usepackage{amssymb,latexsym,amsmath,amsfonts}
\usepackage{amsmath,amscd,amsthm,amssymb}

\numberwithin{equation}{section}

\def\whitebox{{\hbox{\hskip 1pt
 \vrule height 6pt depth 1.5pt
 \lower 1.5pt\vbox to 7.5pt{\hrule width
    3.2pt\vfill\hrule width 3.2pt}%
 \vrule height 6pt depth 1.5pt
 \hskip 1pt } }}
\def\qed{\ifhmode\allowbreak\else\nobreak\fi\hfill\quad\nobreak
     \whitebox\medbreak}
\newcommand{\pf}{\noindent{\bf Proof:}\ }

\newcommand{\ignore}[1]{}

\newcommand{\eproof}{\hfill$\Box$\vspace{4mm}}

\newcommand{\ga}{\gamma}
\newcommand{\F}{\mathbb{F}}
\newcommand{\Z}{\mathbb{Z}}

\newcommand{\Q}{\mathbb{Q}}
\newcommand{\la}{\langle}
\newcommand{\ra}{\rangle}
\newcommand{\eps}{\epsilon}
\newcommand{\Cc}{{\mathbb C}}

\newcommand{\Trace}{{\rm Tr}}

\newtheorem{thm}{Theorem}[section]

\newtheorem{corollary}[thm]{Corollary}

\newtheorem{example}[thm]{Example}
\numberwithin{equation}{section}

\begin{document}

\title[Strongly Regular Cayley Graphs From Cyclotomy]
{Strongly Regular Graphs From Unions of Cyclotomic Classes}

\author[Feng and Xiang]{Tao Feng$^*$, Qing Xiang$^\dagger$}

\thanks{$^*$Supported in part by the Fundamental Research Funds for the central universities.}

\thanks{$^\dagger$Research supported in part by NSF Grant DMS 1001557, and by the Overseas Cooperation Fund (grant 10928101) of China.}

\address{Department of Mathematical Sciences, University of Delaware, Newark, DE 19716, USA.
{\bf Current Address:} Department of Mathematics, Zhejiang University, Hangzhou 310027, Zhejiang, China}
\email{pku.tfeng@yahoo.com.cn}

\address{Department of Mathematical Sciences, University of Delaware, Newark, DE 19716, USA} \email{xiang@math.udel.edu}

\keywords{Cyclotomy, Gauss sum, index 2 Gauss sum, strongly regular graph, uniform cyclotomy.}

\begin{abstract} We give two constructions of strongly regular Cayley graphs on finite fields $\F_q$ by using union of cyclotomic classes and index $2$ Gauss sums. In particular, we obtain twelve infinite families of strongly regular graphs with new parameters.  
\end{abstract}

\maketitle

\section{Introduction}
All graphs considered in this paper are simple and undirected. A {\em strongly regular graph srg} $(v,k,\lambda,\mu)$ is a graph with $v$ vertices that is not complete or edgeless and that has the following properties:
\begin{enumerate}
\item Each vertex is adjacent to $k$ vertices, i.e., the graph is regular of valency $k$,
\item For any two adjacent vertices $x,y$, there are exactly $\lambda$ vertices adjacent to both $x$ and $y$.
\item For any two nonadjacent vertices $x,y$, there are exactly $\mu$ vertices adjacent to both $x$ and $y$.
\end{enumerate}
Classical examples of strongly regular graphs include the Paley graphs. Let $q=4t+1$ be a prime power. The {\it Paley graph} Paley($q$) is the graph with the finite field $\F_q$ as vertex set, where two vertices are adjacent when they differ by a nonzero square. One can check that Paley($q$) is an srg $(4t+1,2t,t-1,t)$. For a survey on strongly regular graphs, we refer the reader to the lecture notes by Brouwer and Haemers \cite{bh}; see also \cite{cg}. Strongly regular graphs are also closely related to two-weight codes, two-intersection sets in finite geometry, and partial difference sets. For these connections, we refer the reader to \cite{CK86, Ma94}. Let $\Gamma$ be a graph. The adjacency matrix of $\Gamma$ is the $(0,1)$-matrix $A$ indexed by the vertex set $V\Gamma$ of $\Gamma$, where $A_{xy} = 1$ when there is an edge between $x$ and $y$ in $\Gamma$ and $A_{xy} = 0$ otherwise. A useful way to check whether a graph is strongly regular is by using the eigenvalues of its adjacency matrix. For convenience we call an eigenvalue {\it restricted} if it has an eigenvector perpendicular to the all-ones vector ${\bf 1}$. 

\begin{thm}\label{char}
For a graph $\Gamma$ of order $v$, not complete or edgeless, with adjacency matrix $A$, the following are equivalent:
\begin{enumerate}
\item $\Gamma$ is an srg $(v, k, \lambda, \mu)$ for certain integers $k, \lambda, \mu$,
\item $A^2 =(\lambda-\mu)A+(k-\mu) I+\mu J$ for certain real numbers $k,\lambda, \mu$, where $I, J$ are the identity matrix and the all-ones matrix, respectively, 
\item $A$ has precisely two distinct restricted eigenvalues.
\end{enumerate}
\end{thm}

This is Theorem~9.1.2 in \cite[p.~115]{bh}. The Paley graphs are probably the simplest examples of the so-called cyclotomic strongly regular graphs, which we define below. Let $p$ be a prime, and $f$ be a positive integer. Let $\F_{p^f}$ be the finite field of order $p^f$, $D\subset \F_{p^f}$ be such that $-D=D$ and $0\not\in D$. We consider the graph $\Gamma$ with the elements of $\F_{p^f}$ as vertices; two vertices are adjacent if and only if their difference belongs to $D$.  That is, $\Gamma$ is the Cayley graph on $\F_{p^f}$ with ``connection" set $D$, usually written as ${\rm Cay}(\F_{p^f}, D)$. When $D$ is a subgroup of the multiplicative group $\F_{p^f}^*$ of $\F_{p^f}$, and if $\Gamma={\rm Cay}(\F_{p^f}, D)$ is strongly regular, then we speak of a {\it cyclotomic strongly regular graph}. The problem of classifying all cyclotomic strongly regular graphs is a venerable one.  We refer the reader to \cite{McE, Lang1, schwh} for detailed studies of this problem. In particular, Schmidt and White \cite{schwh} proposed a conjectural classification of all cyclotomic srg. In this paper, we are interested in constructing srg ${\rm Cay}(\F_{p^f}, D)$, in which $D$ is a union of at least two cosets of a subgroup of $\F_{p^f}^*$ (while a single coset does not give rise to an srg). There are some known examples of such srg. To describe these examples,  we fix a primitive element $\gamma$ of $\F_{p^f}$; let $N$ be an integer greater than 1 such that $N|(p^f-1)$, $C_0$ be the subgroup of $\F_{p^f}^*$ of index $N$, and $C_i=\gamma^i C_0$, $1\leq i\leq N-1$.

\begin{example} {\em (De Lange, \cite{del})} \label{dela}
Let $p=2$, $f=12$, $N=45$. Then 
${\rm Cay}(\F_{p^f}, C_0\cup C_5\cup C_{10})$
is an srg, while ${\rm Cay}(\F_{p^f},C_0)$ is not.
\end{example}

\begin{example}{\em  (Ikuta and Munemasa, \cite{ikutam})}\label{ikutam1} 
Let $p=2$, $f=20$, $N=75$. Then
${\rm Cay}(\F_{p^f}, C_0\cup C_3\cup C_{6}\cup C_9\cup C_{12})$
is an srg, while ${\rm Cay}(\F_{p^f},C_0)$ is not.
\end{example}

\begin{example} {\em  (Ikuta and Munemasa, \cite{ikutam})}\label{ikutam2} 
Let $p=2$, $f=21$, $N=49$. Then
${\rm Cay}(\F_{p^f}, C_0\cup C_1\cup C_{2}\cup C_3\cup C_{4}\cup C_5\cup C_6)$
is an srg, while ${\rm Cay}(\F_{p^f},C_0)$ is not.
\end{example}

We will generalize each of the above three examples into an infinite family. Moreover we obtain nine more infinite families of srg with new parameters by using union of cyclotomic classes.

\section{Gauss sums}
Let $p$ be a prime, $f$ a positive integer, and $q=p^f$. Let $\xi_{p}$ be a fixed complex
primitive $p$th root of unity and let $\Trace_{q/p}$ be the trace from
$\F_{q}$ to $\Z/p\Z$. Define
\begin{equation}\label{defaddchar}
\psi: \F_{q} \rightarrow \Cc^{*}, \quad \psi(x)=\xi_{p}^{\Trace_{q/p}(x)},
\end{equation}
which is easily seen to be a nontrivial character of the additive
group of $\F_{q}$. Let
$$\chi:\F_{q}^* \rightarrow \Cc^{*} $$
be a character of $\F_{q}^{*}$. We define the {\it Gauss sum} by
$$ g(\chi)=\sum_{a \in \F_{q}^*} \chi(a)\psi(a).$$
Note that if $\chi_0$ is the trivial multiplicative character of
$\F_q$, then $g(\chi_0)=-1$. We are usually concerned with nontrivial Gauss sums $g(\chi)$, i.e., those with $\chi\neq \chi_0$. Gauss sums can be viewed as the
Fourier coefficients in the Fourier expansion of $\psi|_{\F_q^*}$
in terms of the multiplicative characters of
 $\F_q$. That is, for every $c\in \F_q^*$,
\begin{equation}\label{inv}
\psi(c)=\frac{1} {q-1}\sum_{\chi\in {\hat \F_q^*}}g({\bar \chi})\chi(c),
\end{equation}
where ${\bar \chi}=\chi^{-1}$ and ${\widehat \F_q^*}$ denotes the character group of $\F_q^*$.

While it is easy to show that the absolute value of a nontrivial Gauss sum $g(\chi)$ is equal to $\sqrt{q}$, the explicit determination of Gauss sums is a difficult problem. However, there are a few cases where the Gauss sums $g(\chi)$ can be explicitly evaluated.  The simplest case is the so-called {\it semi-primitive case}, where there exists an integer $j$ such that $p^j\equiv -1$ (mod $N$) and $N$ is the order of $\chi$ in ${\widehat \F_q^*}$. Some authors \cite{BMW, bew} also refer to this case as uniform cyclotomy, or pure Gauss sums. For future use, we state the following theorem dealing with the semi-primitive case.

\begin{thm} {\em (\cite[p.~364]{bew})}\label{semi}
Let $p$ be a prime, and $N>2$ be an integer. Suppose that there is a positive integer $t$ such that $p^t\equiv -1$ (mod $N$), with $t$ chosen minimal. Let $\chi$ be a character of order $N$ of $\F_{p^r}^*$ for some positive integer $r$. Then $r=2ts$ for some positive integer $s$, and 
\[p^{-r/2}g(\chi)=\begin{cases}(-1)^{s-1},\quad\quad\quad\quad \textup{ if }p=2,\\(-1)^{s-1+\frac{(p^t+1)s}{N}},\quad \textup{ if } p>2.\end{cases}\]
\end{thm}

The next interesting case is the index 2 case, where $-1$ is not in the subgroup $\la p\ra$, the cyclic group generated by $p$ (the characteristic of the finite field $\F_q$),  and $\la p\ra$ has index 2 in $(\Z/N\Z)^*$ (again here $N$ is the order of $\chi$ in ${\widehat \F_q^*}$). Many authors have studied this case, including McEliece \cite{McE}, Langevin \cite{Lang}, Mbodj \cite{Mbo}, Meijer and Van der Vlugt \cite{mv}, and Yang and Xia  \cite{yx}. We state here some results in the index 2 case which will be used in our constructions of strongly regular graphs. Below we use $\phi(N)$ to denote the number of integers $k$ with $1\leq k\leq N$ such that $\gcd(k,N)=1$, and ${\rm ord}_N(p)$ to denote the order of $p$ modulo $N$, which is the smallest positive integer $f$ such that $p^f\equiv 1$ (mod $N$).

\begin{thm}\label{lang}{\em (Langevin \cite{Lang})}
Let $N=p_1^m$, where $m$ is a positive integer, $p_1$ is a prime such that $p_1>3$ and $p_1\equiv 3$ (mod 4). Let $p$ be a prime such that $[(\Z/N\Z)^*:\la p\ra]=2$ (that is, $f:={\rm ord}_N(p)=\phi(N)/2$) and let $q=p^f$. Let $\chi$ be a  multiplicative character of order $N$ of $\F_q$, and $h$ be the class number of $\Q(\sqrt{-p_1})$. Then the Gauss sum $g(\chi)$ over $\F_q$ is determined up to complex conjugation by
$$g(\chi)=\frac{b+c\sqrt{-p_1}}{2}p^{h_0},$$
where 
\begin{enumerate}
\item $h_0=\frac {f-h}{2}$,
\item $b,c\not\equiv 0$ (mod $p$),
\item $b^2+p_1c^2=4p^{h}$,
\item $bp^{h_0}\equiv -2$ (mod $p_1$).
\end{enumerate}
\end{thm}

\begin{thm}\label{mbodj}{\em (Mbodj \cite{Mbo})}
Let $N=p_1^mp_2^{n}$, where $m,n$ are positive integers, $p_1$ and $p_2$ are primes such that $\{p_1\; ({\rm mod}\; 4), p_2\; ({\rm mod}\; 4)\}=\{1,3\}$, ${\rm ord}_{p_1^m}(p)=\phi(p_1^m)$, ${\rm ord}_{p_2^{n}}(p)=\phi(p_2^{n})$. Let $p$ be a prime such that $[(\Z/N\Z)^*:\la p\ra]=2$ (that is, $f:={\rm ord}_N(p)=\phi(N)/2$) and let $q=p^f$. Let $\chi$ be a multiplicative character of order $N$ of $\F_q$, and $h$ be the class number of $\Q(\sqrt{-p_1p_2})$. Then the Gauss sum $g(\chi)$ over $\F_q$ is determined up to complex conjugation by
$$g(\chi)=\frac{b+c\sqrt{-p_1p_2}}{2}p^{h_0},$$
where 
\begin{enumerate}
\item $h_0=\frac {f-h}{2}$,
\item $b,c\not\equiv 0$ (mod $p$),
\item $b^2+p_1p_2c^2=4p^{h}$,
\item $b\equiv 2p^{h/2}$ (mod $\ell$), here $\ell\in\{p_1,p_2\}$ is the prime congruent to 3 modulo 4.
\end{enumerate}
\end{thm}

\section{Cyclotomic classes and periods}
Let $q=p^f$ be a prime power, and let $\ga$ be a fixed primitive element of $\F_q$. Let $N>1$ be a divisor of $q-1$. We define the $N$-th {\em cyclotomic classes}
$C_0,C_1,\ldots, C_{N-1}$ by
$$C_i=\{\ga^{jN+i}\mid 0\leq j\leq \frac {q-1}{N}-1\},$$ where $0\leq i\leq N-1$. That is, $C_0$ is the subgroup of $\F_q^*$ consisting of all nonzero $N$-th powers in $\F_q$, and $C_i=\gamma^i C_0$, $1\leq i\leq N-1$.

Let $\psi$ be the additive character of $\F_q$ defined in (\ref{defaddchar}). The $N$-th {\it cyclotomic periods} are defined by
$$\eta_a=\sum_{x\in C_a}\psi(x),$$
 where $0\leq a\leq N-1$. The relationship between Gauss sums and cyclotomic periods are given as follows. Using (\ref{inv}), we have
\begin{eqnarray*}\label{combin}
\eta_a  & = &\sum_{x \in C_0}\psi(\gamma^a x)\\
& = &\sum_{x \in C_0} \frac{1}{q-1}\sum_{\chi \in \widehat{\mathbb{F}_q^*} } g(\bar{\chi}) \chi(\gamma^a x)\\
& =&\frac{1}{(q-1)}\sum_{\chi \in \widehat{\mathbb{F}_q^*} } g(\bar{\chi})\chi(\gamma^a)\sum_{x \in C_0} \chi(x)\\
& =&\frac{1}{N}\sum_{\chi \in C_0^\perp }g(\bar{\chi})\chi(\gamma^a)
\end{eqnarray*}
where $C_0^{\perp}$ is the subgroup of $\widehat{\F_q^*}$ consisting of all $\chi$ which are trivial on $C_0$, i.e., $C_0^{\perp}$ is the unique subgroup of order $N$. This shows that cyclotomic periods (multiplied by $N$) are linear combinations of Gauss sums with coefficients being 
complex $N$-th roots of unity.

As we already mentioned in Section 2,  the case when ${\rm ord}_N(p)=\phi(N)/2$ and $-1\not\in\la p\ra\leq \Z_N^*$ is usually called {\it the index 2 case}. It is an easy exercise to show that in the index 2 case, $N$ has at most two odd prime divisors. Assume that $N$ is odd, we have the following three possibilities in the index 2 case (see \cite{yx}). Below both $p_1$ and $p_2$ are prime.
\begin{enumerate}
\item $N=p_1^m$, $p_1\equiv 3$ (mod 4);
\item $N=p_1^mp_2^{n}$, $\{p_1\; ({\rm mod}\; 4), p_2\; ({\rm mod}\; 4)\}=\{1,3\}$, ${\rm ord}_{p_1^m}(p)=\phi(p_1^m)$, ${\rm ord}_{p_2^{n}}(p)=\phi(p_2^{n})$;
\item $N=p_1^mp_2^{n}$, $p_1\equiv 1,3$ (mod 4), ${\rm ord}_{p_1^m}(p)=\phi(p_1^m)$ and $p_2\equiv 3$ (mod 4), ${\rm ord}_{p_2^{n}}(p)=\phi(p_2^{n})/2$.
\end{enumerate}
In the following two sections we will consider the cases where new strongly regular graphs are constructed by taking union of cyclotomic classes.

\section{The index 2 case with $N=p_1^m$}

In this section we assume that $N=p_1^m$  (here $m\geq 1$,  $p_1>3$ is a prime such that $p_1\equiv 3$ (mod 4)), $p$ is a prime such that $\gcd(N,p)=1$, and  $f:={\rm ord}_{N}(p)=\phi(N)/2$. Let $q=p^f$, and as before let $C_0,C_1,\ldots ,C_{N-1}$ be the $N$-th cyclotomic classes of $\F_q$. Note that $-C_i=C_i$ for all $0\leq i\leq N-1$ since either $2N|(q-1)$ or $q$ is even. Define
\begin{equation}\label{consecindices}
D=\bigcup_{i=0}^{p_1^{m-1}-1}C_i
\end{equation}
Using $D$ as connection set, we construct the Cayley graph ${\rm Cay}(\F_q, D)$.

\begin{thm}\label{eigenvalues}
The Cayley graph ${\rm Cay}(\F_q, D)$ is a regular graph of valency $|D|$, and it has at most three distinct restricted eigenvalues.
\end{thm} 

\pf Since $-D=D$ and $0\not\in D$, the Cayley graph ${\rm Cay}(\F_q, D)$ is undirected and without loops. It is also regular of valency $|D|$. The restricted eigenvalues of this Cayley graph, as explained in \cite{bwx} (see also \cite[p.~134]{bh}), are 
$$\psi(\ga^a D):=\sum_{x\in D}\psi(\ga^a x),$$
where $\ga$ is a fixed primitive element of $\F_q$, $\psi$ is the additive character of $\F_q$ defined in  (\ref{defaddchar}) and $0\leq a\leq N-1$.

We have 
\begin{eqnarray*}
\psi(\ga^a D)&=&\sum_{i=0}^{p_1^{m-1}-1}\psi(\ga^a C_i)\\
                    &=&\sum_{i=0}^{p_1^{m-1}-1}\eta_{a+i}\\
                    &=&\frac{1}{N}\sum_{\chi\in C_0^{\perp}}g({\bar \chi})\sum_{i=0}^{p_1^{m-1}-1}\chi(\gamma^{a+i}),
\end{eqnarray*}
where $C_0^{\perp}$ is the unique subgroup of ${\widehat \F_q^*}$ of order $N=p_1^m$. For convenience we define 
\[T_a=\sum_{\chi\in C_0^{\perp}}g({\bar \chi})\sum_{i=0}^{p_1^{m-1}-1}\chi(\ga^{a+i}). \]
If $\chi\in C_0^{\perp}$ and  $o(\chi)=1$, then $g({\bar \chi})=-1$, and $\sum_{i=0}^{p_1^{m-1}-1}\chi(\ga^{a+i})=p_1^{m-1}$. If $\chi\in C_0^{\perp}$ and $o(\chi)=p_1^i\; (1\leq i\leq m-1)$, then  $\sum_{i=0}^{p_1^{m-1}-1}\chi(\ga^{a+i})=\chi(\ga^a)\frac{\chi(\ga)^{p_1^{m-1}}-1}{\chi(\ga)-1}=0$. Therefore we have 
$$T_a=-p_1^{m-1}+\sum_{t\in\Z_{p_1^m}^*}g({\bar \chi}^t)\sum_{i=0}^{p_1^{m-1}-1}\chi^t(\ga^{a+i}),$$ where $\chi$ is the character of $\F_q^*$ defined by
\begin{equation}\label{defmulti}
\chi(\ga)=\text{exp}(\frac{2\pi i}{N}).
\end{equation}
By Theorem~\ref{lang}, we have
\[g({\bar \chi})=\frac{b+c\sqrt{-p_1}}{2}p^{h_0},\;\; b,c\not\equiv 0\pmod{p},\]
where $h_0=\frac {f-h}{2}$ and $h$ is the class number of $\Q(\sqrt{-p_1})$, $b^2+p_1c^2=4p^{h}$, and $bp^{h_0}\equiv -2\pmod{p_1}$.  It follows that for any $t\in \Z_{p_1^m}^*$,  $g({\bar \chi}^t)=\frac{b+c(\frac{t}{p_1})\sqrt{-p_1}}{2}p^{h_0}, $ where $(\frac {.}{p_1})$ is the Lengdre symbol.\\

For each $a$, $0\leq a\leq N-1$, there is a unique $i_a\in\{0,1,\ldots,p_1^{m-1}-1\}$, such that $p_1^{m-1}|(a+i_a)$. Write $a+i_a=p_1^{m-1}j_a$ (e.g., when $N=p_1$, we simply have $j_a=a$). For convenience, we introduce the Kronecker delta $\delta_{j_a}$, which equals $1$ if $p_1|j_a$, and $0$ otherwise.  For each $t\in \Z_{p_1^m}^*$, we write $t=t_1+p_1t_2$, where $t_1\in\Z_{p_1}^*$ and $t_2\in\Z_{p_1^{m-1}}$. Now we can compute
\begin{align*}T_a+p_1^{m-1}=&\sum_{t\in\Z_{p_1^m}^*}g({\bar \chi}^t)\sum_{i=0}^{p_1^{m-1}-1}\chi^t(\ga^{a+i})\\
=&p^{h_0}\sum_{t_1\in\Z_{p_1}^*}\sum_{i=0}^{p_1^{m-1}-1}\chi^{t_1}(\ga^{a+i})\frac{b+c(\frac{t_1}{p_1})\sqrt{-p_1}}{2}\sum_{t_2\in \Z_{p_1^{m-1}}}\chi^{p_1t_2}(\ga^{a+i})\\
=&p^{h_0}p_1^{m-1}\sum_{t_1\in\Z_{p_1}^*}\chi^{p_1^{m-1}t_1}(\ga^{j_a})\frac{b+c(\frac{t_1}{p_1})\sqrt{-p_1}}{2}\\
=&\frac{p^{h_0}p_1^{m-1}b}{2}\sum_{t_1\in\Z_{p_1}^*}\chi^{p_1^{m-1}t_1}(\ga^{j_a})+\frac{p^{h_0}p_1^{m-1}c\sqrt{-p_1}}{2}\sum_{t_1\in\Z_{p_1}^*}\chi^{p_1^{m-1}t_1}(\ga^{j_a})\bigg(\frac{t_1}{p_1}\bigg)\\
=&\frac{p^{h_0}p_1^{m-1}b}{2}(p_1\delta_{j_a}-1)+\frac{p^{h_0}p_1^{m-1}c\sqrt{-p_1}}{2}\bigg(\frac{j_a}{p_1}\bigg)\sqrt{-p_1}\\
=&\frac{p^{h_0}p_1^{m-1}b}{2}(p_1\delta_{j_a}-1)-\frac{p^{h_0}p_1^mc}{2}\bigg(\frac{j_a}{p_1}\bigg).
\end{align*}
We remark that when $N=p_1$ (i.e., $m=1$), the second line in the above computations needs to be deleted; everything else still holds true in this case. Therefore we have
\begin{displaymath}
T_a+p_1^{m-1}=
\begin{cases}
\frac{p^{h_0}p_1^{m}b}{2}-\frac{p^{h_0}p_1^{m-1}b}{2}, &\text{if}\; (\frac{j_a}{p_1})=0,\\
\pm\frac{p^{h_0}p_1^{m}c}{2}-\frac{p^{h_0}p_1^{m-1}b}{2} , &\text{if}\;(\frac{j_a}{p_1})\neq 0.\\
\end{cases}
\end{displaymath}
The eigenvalues of ${\rm Cay}(\F_q,D)$ are $|D|=p_1^{m-1}\frac {q-1}{N}=\frac {p^f-1}{p_1}$, and 
\begin{displaymath}
\psi(\ga^aD)=\frac{1}{N}T_a=
\begin{cases}
\frac{p^{h_0}b}{2}-\frac{p^{h_0}b}{2p_1}-\frac{1}{p_1}, &\text{if}\; (\frac{j_a}{p_1})=0,\\
\pm\frac{p^{h_0}c}{2}-\frac{p^{h_0}b}{2p_1}-\frac{1}{p_1} , &\text{if}\;(\frac{j_a}{p_1})\neq 0,\\
\end{cases}
\end{displaymath}
where $0\leq a\leq N-1$. So ${\rm Cay}(\F_q,D)$ has at most three distinct restricted eigenvalues. The proof is now complete.
\eproof

Let $\chi$ be the multiplicative character defined in (\ref{defmulti}), and let
\begin{equation}\label{evaluation}
g({\bar \chi})=\frac{b+c\sqrt{-p_1}}{2}p^{h_0},\;\; b,c\not\equiv 0\pmod{p},
\end{equation}
where $h_0=\frac {f-h}{2}$ and $h$ is the class number of $\Q(\sqrt{-p_1})$, $b^2+p_1c^2=4p^{h}$, and $bp^{h_0}\equiv -2\pmod{p_1}$. We note that while $c$ can only be determined up to sign, $b$ is uniquely determined (without sign ambiguity) by the condition that $bp^{h_0}\equiv -2\pmod{p_1}$ . We have the following corollary.

\begin{corollary}\label{srg1}
Using the above notation, ${\rm Cay}(\F_q,D)$ is a strongly regular graph if and only if $b,c\in \{1,-1\}$.
\end{corollary}

\pf If ${\rm Cay}(\F_q,D)$ is a strongly regular graph, then by Theorem~\ref{char} it has precisely two distinct restricted eigenvalues, $r$ and $s$. As usual, we use $r$ to denote the positive eigenvalue, and $s$ the negative one. By Theorem~\ref{eigenvalues} and the explicit computations of eigenvalues in its proof, we must have $c=\pm b$. Since gcd$(b,c)$ divides $4p^{h}$ and $b,c\not\equiv 0$ (mod $p$), the condition that $c=\pm b$ is equivalent to $b,c\in\{1,-1\}$ or $b,c\in\{2,-2\}$. It is impossible to have $b,c\in\{2,-2\}$: otherwise, from $1+p_1=p^{h}$ we deduce that $p=2$, contradicting the fact $b,c\not\equiv 0\pmod{p}$. Therefore we conclude that $b,c\in\{1,-1\}$.

Fo the converse, noting that if $b,c\in\{1,-1\}$, then $\psi(\gamma^aD)$, $0\leq a\leq N-1$, take only two distinct values. Hence ${\rm Cay}(\F_q,D)$ has precisely two distinct restricted eigenvalues. By Theorem~\ref{char}, ${\rm Cay}(\F_q,D)$ is strongly regular. The proof is now complete.
\eproof

Now if $N=p_1^m$, $p_1$ is a prime congruent to 3 modulo 4, $p_1>3$, and $\frac{1+p_1}{4}=p^{h}$ for some prime $p$, where $h$ is the class number of $\Q(\sqrt{-p_1})$, and  $f:=\text{ord}_N(p)=\phi(N)/2$, then the only possible $b,c$ satisfying (\ref{evaluation}) must be $\pm 1$. This can be seen as follows: from $b^2+p_1c^2=4p^h=1+p_1$ and $b,c\not\equiv 0\pmod{p}$, we get that $1+p_1\leq b^2+p_1c^2=1+p_1$,  hence $b,c\in\{1,-1\}$. Therefore by Corollary~\ref{srg1}, under the above assumptions, the Cayley graph ${\rm Cay}(\F_q, D)$ is strongly regular. Also we note that under the above assumptions, we can further determine when $b$ is equal to 1 and when $b$ is equal to $-1$: First of all, we note that since $[\Z_{p_1^{m}}^* : \la p\ra]=2$, we have $(\frac{p}{p_1})=1$. When $p_1\equiv 3\pmod{8}$, we have $(\frac{2}{p_1})=-1$;  raising both sides of $p^{h_0}b\equiv -2\pmod{p_1}$ to the $\frac{p_1-1}{2}$th power, we obtain $b(\frac{p}{p_1})^{h_0}=(\frac{-1}{p_1})(\frac{2}{p_1})=1$, from which we get $b=1$. Similarly, when $p_1\equiv 7\pmod{8}$, we get $b=-1$. We now set out to find explicit examples of strongly regular Cayley graphs in this way. In the following we only list the examples with $m\geq 2$ since the $m=1$ case was considered previously by Langevin in \cite{Lang1}; see also \cite{schwh}.

\begin{example}\label{1st}
{\em Let $p=2$, $p_1=7$, $N=p_1^m$, $m\geq 2$ is an integer. It is straightforward to check that ${\rm ord}_{7^2}(2)=21=\phi(7^2)/2$. One can easily prove by induction that ${\rm ord}_{N}(2)=\phi(7^m)/2$ for all $m\geq 2$. The class number $h$ of $\Q(\sqrt{-7})$ is equal to 1 (c.f. \cite[p.~514]{hcohen}). Therefore we indeed have $\frac{1+p_1}{4}=p^{h}$ in this case.  We have $f=3\cdot 7^{m-1}$, $b=-1$, $h_0=\frac{f-1}{2}$.  Therefore we obtain a strongly regular Cayley graph ${\rm Cay}(\F_q, D)$, with $v=q=2^{3\cdot 7^{m-1}}$, $k=\frac{v-1}{7}$, and with restricted eigenvalues $r=\frac{2^{h_0+2}-1}{7}$, $s=\frac{-3\cdot 2^{h_0}-1}{7}$.}
\end{example}

We remark that when $m=2$, the srg in Example~\ref{1st} is the same as ${\rm Cay}(\F_{2^{21}}, R_1)$ in Example 3 of \cite{ikutam}.

\begin{example}\label{2nd}
{\em Let $p=3$, $p_1=107$, $N=p_1^m$, $m\geq 2$ is an integer. It is straightforward to check that ${\rm ord}_{107^2}(3)=5671=\phi(107^2)/2$. One can easily prove by induction that ${\rm ord}_{N}(3)=\phi(107^m)/2$ for all $m\geq 2$. The class number $h$ of $\Q(\sqrt{-107})$ is equal to 3 (c.f. \cite[p.~514]{hcohen}). Therefore we indeed have $\frac{1+p_1}{4}=p^{h}$ in this case. We have $f=53\cdot 107^{m-1}$, $b=1$, $h_0=\frac{f-3}{2}$.  Therefore we obtain a strongly regular Cayley graph ${\rm Cay}(\F_q, D)$, with $v=q=3^{53\cdot 107^{m-1}}$, $k=\frac{v-1}{107}$, and with restricted eigenvalues $r=\frac{53\cdot 3^{h_0}-1}{107}$, $s=\frac{-54\cdot 3^{h_0}-1}{107}$.}
\end{example}

\begin{example}\label{3rd}
{\em Let $p=5$, $p_1=19$, $N=p_1^m$,  $m\geq 2$ is an integer. It is straightforward to check that ${\rm ord}_{19^2}(5)=171=\phi(19^2)/2$. One can easily prove by induction that ${\rm ord}_{N}(5)=\phi(19^m)/2$ for all $m\geq 2$. The class number $h$ of $\Q(\sqrt{-19})$ is equal to 1 (c.f. \cite[p.~514]{hcohen}). Therefore we indeed have $\frac{1+p_1}{4}=p^{h}$ in this case. We have $f=9\cdot 19^{m-1}$, $b=1$, $h_0=\frac{f-1}{2}$. Therefore we obtain a strongly regular Cayley graph ${\rm Cay}(\F_q, D)$, with $v=q=5^{9\cdot 19^{m-1}}$, $k=\frac{v-1}{19}$, and with restricted eigenvalues $r=\frac{9\cdot 5^{h_0}-1}{19}$, $s=\frac{-10\cdot 5^{h_0}-1}{19}$.}
\end{example}

\begin{example}\label{4th}
{\em Let $p=5$, $p_1=499$, $N=p_1^m$, $m\geq 2$ is an integer. It is straightforward to check that ${\rm ord}_{499^2}(5)=124251=\phi(499^2)/2$. One can easily prove by induction that ${\rm ord}_{N}(5)=\phi(499^m)/2$ for all $m\geq 2$. The class number $h$ of $\Q(\sqrt{-499})$ is equal to 3 (c.f. \cite[p.~514]{hcohen}). Therefore we indeed have $\frac{1+p_1}{4}=p^{h}$ in this case. We have $f=249\cdot 499^{m-1}$, $b=1$, $h_0=\frac{f-3}{2}$. Therefore we obtain a strongly regular Cayley graph ${\rm Cay}(\F_q, D)$, with $v=q=5^{249\cdot 499^{m-1}}$, $k=\frac{v-1}{499}$, and with restricted eigenvalues $r=\frac{249\cdot 5^{h_0}-1}{499}$, $s=\frac{-250\cdot 5^{h_0}-1}{499}$.}
\end{example}

\begin{example}\label{5th}
{\em Let $p=17$, $p_1=67$,  $N=p_1^m$, $m\geq 2$ is an integer. It is straightforward to check that ${\rm ord}_{67^2}(17)=2211=\phi(67^2)/2$. One can easily prove by induction that ${\rm ord}_{N}(17)=\phi(67^m)/2$ for all $m\geq 2$. The class number $h$ of $\Q(\sqrt{-67})$ is equal to 1 (c.f. \cite[p.~514]{hcohen}). Therefore we indeed have $\frac{1+p_1}{4}=p^{h}$ in this case. We have $f=33\cdot 67^{m-1}$, $b=1$, $h_0=\frac{f-1}{2}$. Therefore we obtain a strongly regular Cayley graph ${\rm Cay}(\F_q, D)$, with $v=q=17^{33\cdot 67^{m-1}}$, $k=\frac{v-1}{67}$, and with restricted eigenvalues $r=\frac{33\cdot 67^{h_0}-1}{67}$, $s=\frac{-34\cdot 67^{h_0}-1}{67}$.}
\end{example}

\begin{example}\label{6th}
{\em Let $p=41$, $p_1=163$, $N=p_1^m$, $m\geq 2$ is an integer. It is straightforward to check that ${\rm ord}_{163^2}(41)=13203=\phi(163^2)/2$. One can easily prove by induction that ${\rm ord}_{N}(41)=\phi(163^m)/2$ for all $m\geq 2$. The class number $h$ of $\Q(\sqrt{-163})$ is equal to 1 (c.f. \cite[p.~514]{hcohen}). Therefore we indeed have $\frac{1+p_1}{4}=p^{h}$ in this case. We have $f=81\cdot 163^{m-1}$, $b=1$, $h_0=\frac{f-1}{2}$. Therefore we obtain a strongly regular Cayley graph ${\rm Cay}(\F_q, D)$, with $v=q=41^{81\cdot 163^{m-1}}$, $k=\frac{v-1}{163}$, and with restricted eigenvalues $r=\frac{81\cdot 41^{h_0}-1}{163}$, $s=\frac{-82\cdot 41^{h_0}-1}{163}$.}
\end{example}

\section{The index 2 case with $N=p_1^mp_2$}

In this section, we assume that $N=p_1^mp_2$ ($m\geq 1$),  $p_1$, $p_2$ are primes such that $\{p_1\; ({\rm mod}\; 4), p_2\; ({\rm mod}\; 4)\}=\{1,3\}$, $p$ is a prime such that $\gcd(p, N)=1$, ${\rm ord}_{p_1^m}(p)=\phi(p_1^m)$ and ${\rm ord}_{p_2}(p)=\phi(p_2)$, and $f:={\rm ord}_N(p)=\phi(N)/2$. Therefore, we are in Case (2), with $n=1$, as listed in the end of Section 3. Let $q=p^f$, and as before let $\gamma$ be a fixed primitive element of $\F_q$,  $C_0=\langle\gamma^N\rangle, C_1=\gamma C_0,\ldots ,C_{N-1}=\gamma^{N-1}C_0$ be the $N$-th cyclotomic classes of $\F_q$. Note that  we have $-C_i=C_i$ for all $0\leq i\leq N-1$ since either $2N|(q-1)$ or $q$ is even. Define
\begin{equation}
D=\bigcup_{i=0}^{p_1^{m-1}-1}C_{ip_2}
\end{equation}
Using $D$ as connection set, we construct the Cayley graph ${\rm Cay}(\F_q, D)$.

\begin{thm}\label{eigenvaluesnew}
The Cayley graph ${\rm Cay}(\F_q, D)$ is a regular graph of valency $|D|$, and it has at most five distinct restricted eigenvalues.
\end{thm}

\pf Since $-D=D$ and $0\not\in D$, the Cayley graph ${\rm Cay}(\F_q, D)$ is undirected and without loops. It is also regular of valency $|D|$. The restricted eigenvalues of this Cayley graph, as explained in \cite{bwx}, are 
$$\psi(\ga^a D):=\sum_{x\in D}\psi(\ga^a x),$$
where $\ga$ is a fixed primitive element of $\F_q$, $\psi$ is the additive character of $\F_q$ defined in  (\ref{defaddchar}) and $0\leq a\leq N-1$.

We have 
\begin{eqnarray*}
\psi(\ga^a D)&=&\sum_{i=0}^{p_1^{m-1}-1}\psi(\ga^a C_{ip_2})\\
                    &=&\sum_{i=0}^{p_1^{m-1}-1}\eta_{a+ip_2}\\
                    &=&\frac{1}{N}\sum_{\chi\in C_0^{\perp}}g({\bar \chi})\chi(\ga^a)\sum_{i=0}^{p_1^{m-1}-1}\chi^{p_2}(\ga^{i}),
\end{eqnarray*}
where $C_0^{\perp}$ is the unique subgroup of ${\widehat \F_q^*}$ of order $N=p_1^mp_2$. For convenience we define 
\[T_a=\sum_{\chi\in C_0^{\perp}}g({\bar \chi})\chi(\ga^a)\sum_{i=0}^{p_1^{m-1}-1}\chi^{p_2}(\ga^{i}).\]
If $\chi\in C_0^{\perp}$, $\chi^{p_2}\neq 1$ and $\chi^{p_1^{m-1}p_2}=1$, then $\sum_{i=0}^{p_1^{m-1}-1}\chi^{p_2}(\ga^i)=\frac{\chi^{p_2p_1^{m-1}}(\ga)-1}{\chi^{p_2}(\ga)-1}=0$. If $\chi\in C_0^{\perp}$ and $\chi^{p_2}=1$, then $\sum_{i=0}^{p_1^{m-1}-1}{\chi}^{p_2}(\ga^i)=p_1^{m-1}$. Therefore,
\[T_a=p_1^{m-1}(-1+A)+B+C,\]
where 
\begin{align*}&A=\sum_{\chi\in C_0^{\perp}\atop{o(\chi)=p_2}}g({\bar \chi}){\chi}(\ga^a),\quad B=\sum_{\chi\in C_0^{\perp}\atop{o(\chi)=p_1^m}}g({\bar \chi}){\chi}(\ga^a)\sum_{i=0}^{p_1^{m-1}-1}{\chi}^{p_2}(\ga^i),\\ &C=\sum_{\chi\in C_0^{\perp}\atop{o(\chi)=N}}g({\bar \chi}){\chi}(\ga^a)\sum_{i=0}^{p_1^{m-1}-1}{\chi}^{p_2}(\ga^i).
\end{align*}
Below we will compute $A$, $B$, $C$ individually. 

For each $a$, $0\leq a\leq N-1$, there is a unique $i_a\in\{0,1,\ldots,p_1^{m-1}-1\}$, such that $p_1^{m-1}|(a+p_2i_a)$. Write $a+p_2i_a=p_1^{m-1}j_a$. Again we introduce the Kronecker delta $\delta_{j_a,p_1}$, which equals 1 if $p_1|j_a$, 0 otherwise. Also we define $\delta_{a,p_2}$ by setting it equal to 1 if $p_2|a$, 0 otherwise.

Since ${\rm ord}_{p_2}(p)=\phi(p_2)$, we have $p^{\frac{p_2 -1}{2}}\equiv -1$ (mod $p_2$). By Theorem \ref{semi}, we have $g({\bar \chi})=(-1)^{\frac{p_1-1}{2}-1}\sqrt{q}$ for each $\chi$ of order $p_2$.  It follows that $$A=(-1)^{\frac{p_1-1}{2}-1}\sqrt{q}\sum_{o(\chi)=p_2}{\chi}(\ga^a)=(-1)^{\frac{p_1-1}{2}-1}\sqrt{q}(p_2\delta_{a,p_2}-1).$$

Similarly,  we have $g({\bar \chi})=(-1)^{\frac{p_2-1}{2}-1}\sqrt{q}$ for each $\chi$ of order $p_1^m$. Let $\chi_1$ be the character of order $p_1^m$ in ${\widehat \F_q^*}$ defined by $\chi_1(\ga)=\textup{exp}(\frac{2\pi i}{p_1^m})$. We have 
\[B=(-1)^{\frac{p_2-1}{2}-1}\sqrt{q}\sum_{i=0}^{p_1^{m-1}-1}\sum_{t\in \Z_{p_1^m}^*}{\chi_1}^{t(p_2i+a)}(\ga)=(-1)^{\frac{p_2-1}{2}-1}\sqrt{q}p_1^{m-1}(p_1\delta_{j_a,p_1}-1).\]

Let $\chi_1$ be defined as above and $\chi_2$ be the character of order $p_2$ in ${\widehat \F_q^*}$ defined by $\chi_2(\ga)=\textup{exp}(\frac{2\pi i}{p_2})$. By Theorem~\ref{mbodj}, we have 
\begin{equation}\label{evalGauss}
g({\bar \chi_1}{\bar \chi_2})=\frac{b+c\sqrt{-p_1p_2}}{2}p^{h_0},
\end{equation}
where $h_0=\frac {f-h}{2}$ ($h$ is the class number of $\Q(\sqrt{-p_1p_2})$), $b,c\not\equiv 0$ (mod $p$), $b^2+p_1p_2c^2=4p^{h}$, and $b\equiv 2p^{h/2}$ (mod $\ell$), here $\ell\in\{p_1,p_2\}$ is the prime congruent to 3 modulo 4.

Every character in $C_0^{\perp}$ of order $p_1^mp_2$ is of the form $\chi_1^u\chi_2^v$ with $u\in\Z_{p_1^m}^*$, $v\in\Z_{p_2}^*$. Let $\sigma_{u,v}$ be the Galois automorphism of $\Q(\xi_p,\xi_N)$ defined by $\sigma_{u,v}(\xi_p)=\xi_p$, $\sigma_{u,v}(\xi_{p_1^m})=\xi_{p_1^m}^u$, $\sigma_{u,v}(\xi_{p_2})=\xi_{p_2}^v$. Then $$g({\bar \chi_1}^u{\bar \chi}_2^v)=\sigma_{u,v}(g({\bar \chi_1}{\bar \chi_2}))=\frac{b+c(\frac{u}{p_1})(\frac{v}{p_2})\sqrt{-p_1p_2}}{2}p^{h_0}.$$ We are now ready to compute $C$.

\begin{align*}
 C&=p^{h_0}\sum_{u\in\Z_{p_1^m}^*}\sum_{v\in \Z_{p_2}^*}\bigg[\frac{b}{2}+\frac{c}{2}\bigg(\frac{u}{p_1}\bigg)\bigg(\frac{v}{p_2}\bigg)\sqrt{-p_1p_2}\bigg]\sum_{i=0}^{p_1^{m-1}-1}{\chi}_1^u(\ga^{p_2i+a}){\chi}_2^v(\ga^a)\\
&=\frac{b}{2}p^{h_0}\bigg(\sum_{i=0}^{p_1^{m-1}-1}\sum_{u\in\Z_{p_1^m}^*}{\chi}_1^u(\ga^{p_2i+a})\bigg)\bigg(\sum_{v\in\Z_{p_2}^*}{\chi}_2^v(\ga^a)\bigg)\\
&\quad+\frac{c}{2}p^{h_0}\sqrt{-p_1p_2}\bigg(\sum_{i=0}^{p_1^{m-1}-1}\sum_{u\in\Z_{p_1^m}^*}\bigg(\frac{u}{p_1}\bigg){\chi}_1^u(\ga^{p_2i+a})\bigg)\bigg( \sum_{v\in\Z_{p_2}^*}\bigg(\frac{v}{p_2}\bigg){\chi}_2^v(\ga^a)\bigg)\\
&=\frac{b}{2}p^{h_0}p_1^{m-1}(p_1\delta_{j_a,p_1}-1)(p_2\delta_{a,p_2}-1)\\
&\quad+\frac{c}{2}p^{h_0}\sqrt{-p_1p_2}\bigg(\sum_{i=0}^{p_1^{m-1}-1}\sum_{u\in\Z_{p_1^m}^*}\bigg(\frac{u}{p_1}\bigg){\chi}_1^u(\ga^{p_2i+a})\bigg)\bigg( \sum_{v\in\Z_{p_2}^*}\bigg(\frac{v}{p_2}\bigg){\chi}_2^v(\ga^a)\bigg).
\end{align*}
We have $\sum_{v\in\Z_{p_2}^*}\big(\frac{v}{p_2}\big){\chi}_2^v(\ga^a)=\big(\frac{a}{p_2}\big)\sqrt{p_2^*}$, where $p_2^*=(-1)^{\frac {p_2-1}{2}}p_2$, and
\begin{align*}
\sum_{i=0}^{p_1^{m-1}-1}\sum_{u\in\Z_{p_1^m}^*}\bigg(\frac{u}{p_1}\bigg){\chi}_1^u(\ga^{p_2i+a}) 
&=\sum_{i=0}^{p_1^{m-1}-1}\bigg(\sum_{u_1=1}^{p_1-1}\bigg(\frac{u_1}{p_1}\bigg){\chi}_1^{u_1}(\ga^{p_2i+a})\bigg)\bigg(\sum_{u_2\in\Z_{p_1^{m-1}}}{\chi}_1^{p_1u_2}(\ga^{p_2i+a})\bigg)\\
&=p_1^{m-1}\sum_{u_1=1}^{p_1-1}\bigg(\frac{u_1}{p_1}\bigg){\chi}_1^{u_1}(\ga^{p_1^{m-1}j_a})=p_1^{m-1}\bigg(\frac{j_a}{p_1}\bigg)\sqrt{p_1^*},
\end{align*}
where $p_1^*=(-1)^{\frac {p_1-1}{2}}p_1$.
Therefore, putting the above computations together, we have
\begin{align*}
 C&=\frac{b}{2}p^{h_0}p_1^{m-1}(p_1\delta_{j_a,p_1}-1)(p_2\delta_{a,p_2}-1)+\frac{c}{2}p^{h_0}\sqrt{-p_1p_2}\bigg(\frac{a}{p_2}\bigg)\sqrt{p_2^*}\cdot p_1^{m-1}\bigg(\frac{j_a}{p_1}\bigg)\sqrt{p_1^*}\\
&=\frac{b}{2}p^{h_0}p_1^{m-1}(p_1\delta_{j_a,p_1}-1)(p_2\delta_{a,p_2}-1)-\bigg(\frac{a}{p_2}\bigg)\bigg(\frac{j_a}{p_1}\bigg)\frac{c}{2}p^{h_0}p_1^mp_2
\end{align*}
We conclude that
\begin{align*}T_a&=p_1^{m-1}(-1+A)+B+C\\
&=-p_1^{m-1}+p_1^{m-1}(-1)^{\frac{p_1-1}{2}-1}\sqrt{q}(p_2\delta_{a,p_2}-1)+(-1)^{\frac{p_2-1}{2}-1}\sqrt{q}p_1^{m-1}(p_1\delta_{j_a,p_1}-1)\\
&\quad+\frac{b}{2}p^{h_0}p_1^{m-1}(p_1\delta_{j_a,p_1}-1)(p_2\delta_{a,p_2}-1)-\bigg(\frac{a}{p_2}\bigg)\bigg(\frac{j_a}{p_1}\bigg)\frac{c}{2}p^{h_0}p_1^mp_2\\
\end{align*}
Noting that $\{p_1\; ({\rm mod}\; 4), p_2\; ({\rm mod}\; 4)\}=\{1,3\}$, we have
\begin{align*}
&T_a=-p_1^{m-1}-(-1)^{\frac{p_1-1}{2}}p_1^{m-1}p_2\sqrt{q}\delta_{a,p_2}-(-1)^{\frac{p_2-1}{2}}p_1^m\sqrt{q}\delta_{j_a,p_1}\\
&\quad+\frac{b}{2}p^{h_0}p_1^{m-1}(p_1\delta_{j_a,p_1}-1)(p_2\delta_{a,p_2}-1)-\bigg(\frac{a}{p_2}\bigg)\bigg(\frac{j_a}{p_1}\bigg)\frac{c}{2}p^{h_0}p_1^mp_2.
\end{align*}
We consider four cases:
\begin{enumerate}
\item $\delta_{j_a,p_1}=\delta_{a,p_2}=0$. In this case, we have
\begin{align*}T_a&=-p_1^{m-1}+ \frac{b}{2}p^{h_0}p_1^{m-1}-\bigg(\frac{a}{p_2}\bigg)\bigg(\frac{j_a}{p_1}\bigg)\frac{c}{2}p^{h_0}p_1^mp_2.
\end{align*}
 Set
\begin{align*}
  c_+&=-p_1^{m-1}+\frac{b}{2}p^{h_0}p_1^{m-1}+\frac{c}{2}p^{h_0}p_1^mp_2,\\
  c_-&=-p_1^{m-1}+\frac{b}{2}p^{h_0}p_1^{m-1}-\frac{c}{2}p^{h_0}p_1^mp_2.
\end{align*}
Then $T_a=c_+$ or $c_-$ according as $\bigg(\frac{a}{p_2}\bigg)\bigg(\frac{j_a}{p_1}\bigg)=-1$ or $1$.
\item $\delta_{j_a,p_1}=1,\;\delta_{a,p_2}=1$. In this case, we have
\begin{align*}T_a&=-p_1^{m-1}-(-1)^{\frac{p_1-1}{2}}p_1^{m-1}p_2\sqrt{q}-(-1)^{\frac{p_2-1}{2}}p_1^m\sqrt{q}+\frac{b}{2}p^{h_0}p_1^{m-1}(p_1-1)(p_2-1).
\end{align*}
For future use, we will denote this value of $T_a$ by $c_1$.
\item $\delta_{j_a,p_1}=0,\;\delta_{a,p_2}=1$. In this case, we have
\begin{align*}T_a&=-p_1^{m-1}-(-1)^{\frac{p_1-1}{2}}p_1^{m-1}p_2\sqrt{q}-\frac{b}{2}p^{h_0}p_1^{m-1}(p_2-1).
\end{align*}
For future use, we will denote this value of $T_a$ by $c_2$.
\item $\delta_{j_a,p_1}=1,\;\delta_{a,p_2}=0$. In this case, we have
\begin{align*}T_a&=-p_1^{m-1}-(-1)^{\frac{p_2-1}{2}}p_1^m\sqrt{q}-\frac{b}{2}p^{h_0}p_1^{m-1}(p_1-1).
\end{align*}
For future use, we will denote this value of $T_a$ by $c_3$.
\end{enumerate}
 In summary, $T_a$, $0\leq a\leq N-1$, belong to $\{c_+,c_-,c_1,c_2,c_3\}$. Therefore the restricted eigenvalues $\psi(\ga^aD)=T_a/N$ can take at most 5 distinct values. The proof is now complete.
\eproof

We now determine when ${\rm Cay}(\F_q,D)$ is strongly regular. Recall that $h$ is the class number of $\Q(\sqrt{-p_1p_2})$, and $b,c$ appeared in (\ref{evalGauss}).

\begin{corollary}\label{srg2}
Using the above notation, ${\rm Cay}(\F_q,D)$ is a strongly regular graph if and only if $b,c\in \{1,-1\}$, $h$ is even and $p_1=2p^{h/2}+(-1)^{\frac{p_1-1}{2}}b$, $p_2=2p^{h/2}-(-1)^{\frac{p_1-1}{2}}b$. 
\end{corollary}

\pf  If  the graph ${\rm Cay}(\F_q,D)$ is strongly regular, then by Theorem~\ref{char}, it has precisely two distinct restricted eigenvalues. Since $c\not\equiv 0$ (mod $p$), we have $c_+\neq c_-$. Therefore we must have $c_1,c_2,c_3\in\{c_+,c_-\}$. It follows that

\begin{equation}\label{c_1}
-(-1)^{\frac{p_1-1}{2}}2p_2p^{h/2}-(-1)^{\frac{p_2-1}{2}}2p_1p^{h/2}+b(p_1p_2-p_1-p_2)=\eps_1cp_1p_2,
\end{equation}

\begin{equation}\label{c_2}
-(-1)^{\frac{p_1-1}{2}}2p^{h/2}=b+\eps_2cp_1,
\end{equation}

\begin{equation}\label{c_3}
-(-1)^{\frac{p_2-1}{2}}2p^{h/2}=b+\eps_3cp_2,
\end{equation}
for some $\eps_1,\eps_2,\eps_3\in\{1, -1\}$. Hence $h$ must be even. Squaring both sides of (\ref{c_2}),   and recall that $b^2+p_1p_2c^2=4p^h$, we obtain $2b\eps_2+cp_1=cp_2$. Squaring both sides of (\ref{c_3}), we obtain $2b\eps_3+cp_2=cp_1$. Combining these two equations, we have $\eps_3=-\eps_2$.  Now substracting $p_2$ copies of (\ref{c_2}) and $p_1$ copies of (\ref{c_3}) from (\ref{c_1}), we get $b=\eps_1c$. Using the same argument as in the proof of Corollary \ref{srg1} , we deduce that $b,c\in\{1,-1\}$. 

Since $p_1,p_2$ are positive and $b=\pm 1$, from (\ref{c_2}) and (\ref{c_3}) we obtain $-(-1)^{\frac{p_1-1}{2}}=\eps_2c=-\eps_3c$. Consequently, 
\[p_1=2p^{h/2}+(-1)^{\frac{p_1-1}{2}}b, \quad p_2=2p^{h/2}-(-1)^{\frac{p_1-1}{2}}b.\]

For the converse, noting that if $b,c\in\{1,-1\}$, $h$ is even, and $p_1=2p^{h/2}+(-1)^{\frac{p_1-1}{2}}b$, $p_2=2p^{h/2}-(-1)^{\frac{p_1-1}{2}}b$, then with $\eps_1=bc$, $\eps_2=-(-1)^{\frac{p_1-1}{2}}c$, $\eps_3=-\eps_2$, the three equations, (\ref{c_1}), (\ref{c_2}) and (\ref{c_3}), will hold. Therefore $c_1,c_2,c_3\in\{c_+,c_-\}$. It follows that $\psi(\gamma^aD)$, $0\leq a\leq N-1$, take precisely two distinct values $r=\frac{c_+}{N}$, $s=\frac{c_-}{N}$. Hence ${\rm Cay}(\F_q,D)$ has precisely two distinct restricted eigenvalues $r$ and $s$. By Theorem~\ref{char}, ${\rm Cay}(\F_q,D)$ is strongly regular. The proof is now complete.
\eproof

Let us find some concrete examples of srg arising in this way. Recall that the $b$ in (\ref{evalGauss}) is uniquely determined by the condition that $b\equiv 2p^{h/2}$ (mod $\ell$), where $\ell\in\{p_1,p_2\}$ is the prime congruent to 3 modulo 4.

\begin{example}\label{21st}
{\em Let $p=2$, $p_1=3$, $p_2=5$, $N=3^m\cdot 5$, with $m\geq 1$. One can easily prove by induction that $f:={\rm ord}_{N}(2)=\phi(N)/2=4\cdot 3^{m-1}$ for all $m\geq 1$. The class number $h$ of $\Q(\sqrt{-15})$ is equal to 2 (c.f. \cite[p.~514]{hcohen}). Since $1+p_1p_2=4p^h$, we have $b,c\in\{1,-1\}$.  From $b\equiv 2p^{h/2}$ (mod $\ell$), here $\ell=3$, we get $b=1$. The conditions in Corollary \ref{srg2} are all satisfied.  Therefore we obtain a strongly regular Cayley graph ${\rm Cay}(\F_q, D)$, with $$v=q=2^{4\cdot 3^{m-1}},\; k=\frac{v-1}{15}=16^{3^{m-1}-1}+16^{3^{m-1}-2}+\cdots +16+1,$$ and with restricted eigenvalues $r=\frac{2^{h_0+3}-1}{15}$, $s=\frac{-7\cdot 2^{h_0}-1}{15}$, where $h_0=\frac{f-h}{2}=2\cdot 3^{m-1}-1$.}
\end{example}

We remark that when $m=2$, the srg in Example~\ref{21st} is the same as Example~\ref{dela} by De Lange.

\begin{example}\label{22st}
{\em Let $p=2$, $p_1=5$, $p_2=3$, $N=5^m\cdot 3$, with $m\geq 1$. One can easily prove by induction that $f:={\rm ord}_{N}(2)=\phi(N)/2=4\cdot 5^{m-1}$ for all $m\geq 1$. The class number $h$ of $\Q(\sqrt{-15})$ is equal to 2 . Since $1+p_1p_2=4p^h$, we have $b,c\in\{1,-1\}$.  For the same reason as in Example~\ref{21st} we have $b=1$. The conditions in Corollary \ref{srg2} are all satisfied.  Therefore we obtain a strongly regular Cayley graph ${\rm Cay}(\F_q, D)$, with $$v=q=2^{4\cdot 5^{m-1}},\; k=\frac{v-1}{15}=16^{5^{m-1}-1}+16^{5^{m-1}-2}+\cdots +16+1,$$ and with restricted eigenvalues $r=\frac{2^{h_0+3}-1}{15}$, $s=\frac{-7\cdot 2^{h_0}-1}{15}$, where $h_0=\frac{f-h}{2}=2\cdot 5^{m-1}-1$.}
\end{example}

We remark that when $m=2$, the srg in Example~\ref{22st} is the same as Example~\ref{ikutam1} by Ikuta and Munemasa.

\begin{example}\label{23st}
{\em Let $p=3$, $p_1=5$, $p_2=7$, $N=5^m\cdot 7$, with $m\geq 1$. One can easily prove by induction that $f:={\rm ord}_{N}(3)=\phi(N)/2=12\cdot 5^{m-1}$ for all $m\geq 1$. The class number $h$ of $\Q(\sqrt{-35})$ is equal to 2 (c.f. \cite[p.~514]{hcohen}). Since $1+p_1p_2=4p^h$, we have $b,c\in\{1,-1\}$.  From $b\equiv 2p^{h/2}$ (mod $\ell$), here $\ell=7$, we get $b=-1$. The conditions in Corollary \ref{srg2} are all satisfied.  Therefore we obtain a strongly regular Cayley graph ${\rm Cay}(\F_q, D)$, with $v=q=3^{12\cdot 5^{m-1}},\; k=\frac{v-1}{35}$ and with restricted eigenvalues $r=\frac{17\cdot 3^{h_0}-1}{35}$, $s=\frac{-18\cdot 3^{h_0}-1}{35}$, where $h_0=\frac{f-h}{2}=6\cdot 5^{m-1}-1$.}
\end{example}

\begin{example}\label{24st}
{\em Let $p=3$, $p_1=7$, $p_2=5$, $N=7^m\cdot 5$, with $m\geq 1$. One can easily prove by induction that $f:={\rm ord}_{N}(3)=\phi(N)/2=12\cdot 7^{m-1}$ for all $m\geq 1$. The class number $h$ of $\Q(\sqrt{-35})$ is equal to 2 (c.f. \cite[p.~514]{hcohen}). Since $1+p_1p_2=4p^h$, we have $b,c\in\{1,-1\}$.  From the same reason as in Example~\ref{23st} we have $b=-1$. The conditions in Corollary \ref{srg2} are all satisfied.  Therefore we obtain a strongly regular Cayley graph ${\rm Cay}(\F_q, D)$, with $v=q=3^{12\cdot 7^{m-1}},\; k=\frac{v-1}{35}$ and with restricted eigenvalues $r=\frac{17\cdot 3^{h_0}-1}{35}$, $s=\frac{-18\cdot 3^{h_0}-1}{35}$, where $h_0=\frac{f-h}{2}=6\cdot 7^{m-1}-1$.}
\end{example}

\begin{example}\label{25st}
{\em Let $p=3$, $p_1=17$, $p_2=19$, $N=17^m\cdot 19$, with $m\geq 1$. One can easily prove by induction that $f:={\rm ord}_{N}(3)=\phi(N)/2=144\cdot 17^{m-1}$ for all $m\geq 1$. The class number $h$ of $\Q(\sqrt{-323})$ is equal to 4 (c.f. \cite[p.~514]{hcohen}). Since $1+p_1p_2=4p^h$, we have $b,c\in\{1,-1\}$.  From $b\equiv 2p^{h/2}$ (mod $\ell$), here $\ell=19$, we get $b=-1$. The conditions in Corollary \ref{srg2} are all satisfied.  Therefore we obtain a strongly regular Cayley graph ${\rm Cay}(\F_q, D)$, with $v=q=3^{144\cdot 17^{m-1}},\; k=\frac{v-1}{323},$ and with restricted eigenvalues $r=\frac{161\cdot 3^{h_0}-1}{323}$, $s=\frac{-162\cdot 3^{h_0}-1}{323}$, where $h_0=\frac{f-h}{2}=72\cdot 17^{m-1}-2$.}
\end{example}

\begin{example}\label{26st}
{\em Let $p=3$, $p_1=19$, $p_2=17$, $N=19^m\cdot 17$, with $m\geq 1$. One can easily prove by induction that $f:={\rm ord}_{N}(3)=\phi(N)/2=144\cdot 19^{m-1}$ for all $m\geq 1$. The class number $h$ of $\Q(\sqrt{-323})$ is equal to 4 (c.f. \cite[p.~514]{hcohen}). Since $1+p_1p_2=4p^h$, we have $b,c\in\{1,-1\}$.  For the same reason as in Example~\ref{25st} we have $b=-1$. The conditions in Corollary \ref{srg2} are all satisfied.  Therefore we obtain a strongly regular Cayley graph ${\rm Cay}(\F_q, D)$, with $v=q=3^{144\cdot 19^{m-1}},\; k=\frac{v-1}{323},$ and with restricted eigenvalues $r=\frac{161\cdot 3^{h_0}-1}{323}$, $s=\frac{-162\cdot 3^{h_0}-1}{323}$, where $h_0=\frac{f-h}{2}=72\cdot 19^{m-1}-2$.}
\end{example}

\section{Concluding Remarks}

We have constructed strongly regular Cayley graphs on $\F_q$ by using union of cyclotomic classes of $\F_q$ and index 2 Gauss sums. Twelve infinite families of srg with new parameters are obtained in this way. It is natural to ask whether further examples can be found by using Corollary 4.2 and 5.2. One can certainly use a computer to search for more prime pairs $(p,p_1)$ satisfying the conditions of Corollary 4.2, and prime triples $(p,p_1,p_2)$ satisfying the conditions of Corollary 5.2. But we suspect that it is unlikely one can find new examples in view of the computer search performed by White and Schmidt \cite{schwh} and the theoretic results therein. 

Another natural question is whether we get interesting fusion schemes of the cyclotomic association schemes by using the srg arising from Corollary 4.2 and 5.2.  We use the construction in Section 4 to explain this problem in some detail below. 

Let $q=p^f$, where $p$ is a prime and $f$ a positive integer. Let $\gamma$ be a fixed primitive element of $\F_q$ and $N|(q-1)$ with $N>1$. As usual, let $C_0=\langle \gamma^N\rangle$, and $C_i=\gamma^i C_0$, $1\leq i\leq N-1$, be the $N$-th cyclotomic classes of $\F_q$. Assume that $-1\in C_0$. Define $R_0=\{(x,x) \mid  x\in \F_q\}$, and for $i\in \{1,2,\ldots ,N\}$, define $R_i=\{(x,y)\mid x,y\in \F_q, x-y\in C_{i-1}\}$. Then $(\F_q, \{R_i\}_{0\leq i\leq N})$ is a symmetric association scheme, which is called {\it the cyclotomic association scheme of class $N$ over $\F_q$}. Now assume that we are in the situation of Section 4. That is, $N=p_1^m$, where $m\geq 1$,  $p_1>3$ is a prime such that $p_1\equiv 3$ (mod 4); $p$ is a prime such that $\gcd(N,p)=1$, and  $f:={\rm ord}_{N}(p)=\phi(N)/2$. For $0\leq k\leq p_1-1$, define
\begin{equation}
D_k=\bigcup_{i=0}^{p_1^{m-1}-1}C_{i+kp_1^{m-1}}
\end{equation}
Note that $D_0$ is the same as $D$ in (\ref{consecindices}), $D_k=\gamma^{kp_1^{m-1}}D_0$, and $D_0, D_1, \ldots ,D_{p_1-1}$ form a partition of $\F_q^*$. Define $R'_0=R_0$ and
\begin{equation}\label{defREL2}
R'_k=\{(x,y)\mid x,y\in \F_q, x-y\in D_{k-1}\}.
\end{equation} 
It is natural to ask whether $(\F_q, \{R'_k\}_{0\leq k\leq p_1})$ is an association scheme. We give an affirmative answer to this question in a subsequent paper \cite{fengwux}. Also included in \cite{fengwux} are some interesting properties of this fusion scheme in relation to A.V. Ivanov's conjecture \cite{IP, ikutam}.

\end{document}